\documentclass[11pt]{article}
\usepackage{amsmath, amssymb, latexsym, times, scalerel}
\usepackage{amsthm}
\usepackage{rsfso}

\usepackage{amscd}
\usepackage{booktabs} 
\usepackage{wrapfig}
\usepackage{tikz}
\usepackage{pgfplots}
\usepackage{float}
\pgfplotsset{compat=newest}
\usepackage{authblk}
\usepackage{mathrsfs}
\usepackage{bbm}
\usepackage{graphicx}
\usepackage{stmaryrd}
\usepackage[all]{xy}
\usepackage{enumitem}
\usepackage{tikz-cd}

\usepackage{tocloft}
\let\OLDthebibliography\thebibliography
\renewcommand\thebibliography[1]{
  \OLDthebibliography{#1}
  \setlength{\parskip}{4pt}
  \setlength{\itemsep}{2pt plus 0.3ex}
}

\usepackage{hyperref}
\hypersetup{
    colorlinks=true,
    linkcolor=blue,
    filecolor=magenta,      
    urlcolor=cyan,
    citecolor=magenta
}

\usepackage{setspace}
\newtheorem{innercustomgeneric}{\customgenericname}
\providecommand{\customgenericname}{}
\newcommand{\newcustomtheorem}[2]{%
  \newenvironment{#1}[1]
  {%
   \renewcommand\customgenericname{#2}%
   \renewcommand\theinnercustomgeneric{##1}%
   \innercustomgeneric
  }
  {\endinnercustomgeneric}
}

\newcustomtheorem{customthm}{Theorem}
\newcustomtheorem{customconj}{Conjecture}  
\newcustomtheorem{customprop}{Proposition}  
\newcustomtheorem{customcor}{Corollary}

\theoremstyle{plain}
\newtheorem{thm}{Theorem}
\newtheorem{lem}[thm]{Lemma}

\newtheorem{thm-def}[thm]{Theorem-Definition}

\theoremstyle{definition}

\newtheorem{cor-def}[thm]{Corollary-Definition}

\theoremstyle{remark}

\makeatletter
\renewcommand{\paragraph}{%
  \@startsection{paragraph}{4}%
  {\z@}{1.25ex \@plus 1ex \@minus .2ex}{-1em}%
  {\normalfont\normalsize\bfseries}%
}
\makeatother

\usepackage{adjustbox}

\newcommand{\bigboxtimes}{\mathord{\adjustbox{width={0.5cm}, valign=c}{${\boxtimes}$}}}

\newcommand{\gl}{\mathrm{GL}}
\newcommand{\ggl}{\mathfrak{gl}}

\newcommand{\bbv}{\mathbb{V}}

\newcommand{\cf}{\mathcal{F}}
\newcommand{\co}{\mathcal{O}}
\newcommand{\cp}{\mathcal{P}}
\newcommand{\cl}{\mathcal{L}}

\newcommand{\ch}{\mathcal{H}}

\newcommand{\cs}{\mathcal{S}}

\newcommand{\xx}{\mathcal{X}}

\newcommand{\kg}{\mathfrak{g}}

\newcommand{\bn}{\mathbf{N}}
\newcommand{\bz}{\mathbf{Z}}

\newcommand{\bp}{\mathbf{P}}

\newcommand{\rtt}{\mathrm{T}}

\newcommand{\rS}{\mathrm{S}}

\newcommand{\val}{\mathrm{val}}

\newcommand{\Ad}{\mathrm{Ad}}

\newcommand{\diag}{\mathrm{diag}}

\newcommand{\pic}{\mathrm{Pic}}

\newcommand{\spf}{\mathrm{Spf}}

\newcommand{\varep}{\varepsilon}

\newcommand{\sm}{\mathrm{sm}}

\newcommand{\ql}{\mathbf{Q}_{\ell}}

\newcommand{\fq}{\mathbf{F}_{q}}

\newcommand{\ct}{\mathcal{T}}

\newcommand{\cc}{\mathcal{C}}

\def\cb{{\mathcal B}}

\def\kg{{\mathfrak{g}}}

\usepackage{xcolor}
\usepackage{version}

\excludeversion{NB}

\newcommand{\bsm}{\begin{smallmatrix}}
\newcommand{\esm}{\end{smallmatrix}}

\def\<{\langle\,}
\def\>{\,\rangle}

\def\kg{\mathfrak g}

\def\<{\langle}
\def\>{\rangle}
\def\={\equiv}

\def\fq{\mathbb{F}_{q}}

\def\val{{\rm val}}

\def\Ker{\mathrm{Ker\,}}

\author{Zongbin \textsc{Chen}}

\title{\bf On the geometric fundamental lemma of Kottwitz}

\begin{document}
\maketitle

\begin{abstract}

We give a proof of the geometric fundamental lemma of Kottwitz. As explained by Laumon, this implies the fundamental lemma for the unitary groups.

\end{abstract}

\section{The main theorem}

Let $k$ be the finite field $\fq$ of characteristic $p$, {$p$ sufficiently large}. Let $F=k(\!(\varep)\!)$ be the field of Laurent series, $\co=k[\![\varep]\!]$ its ring of integers.
Let $G=\gl_{d}$ and let $\kg=\ggl_{d}$ be its Lie algebra. 
Fix a partition $d=d_{1}+\cdots+d_{r}, d_{i}\in \bn, i=1,\cdots,r$. 
Let $\gamma=\diag(\gamma_{i})_{i=1}^{r}\in \kg(\co)$ be a regular semisimple element, with $\gamma_{i}\in \ggl_{d_{i}}(\co)$ being anisotropic, i.e. the characteristic polynomial $P_{i}(x)$ of $\gamma_{i}$ over $F$ is irreducible.

Recall that the \emph{affine Springer fiber} at $\gamma$ is the closed sub ind-$k$-scheme of the \emph{affine grassmannian} $\xx=G(\!(\varep)\!)/G[\![\varep]\!]$ defined by
$$
\xx_{\gamma}=\big\{[g]\in \xx\,\big|\,\Ad(g^{-1})\gamma \in \kg[\![\varep]\!]\big\}.
$$ 
It is finite dimensional and locally of finite type. 
Let $\rtt$ be the centralizer of $\gamma$ in $G_{F}$, then $\rtt$ acts on the affine Springer fiber $\xx_{\gamma}$ by left translation. 
Let $\rS$ be the maximal $F$-split subtorus of $\rtt$.
Let $\Lambda\subset \rS(F)$ be the subgroup generated by $\chi(\varep),\,\chi\in X_{*}(\rS)$, then $\Lambda$ acts simply transitively on the irreducible components of $\xx_{\gamma}$. 
The connected components of $\xx_{\gamma}$ is naturally parametrized by $\bz$ with the morphism
$$
\xx_{\gamma}\to \bz,\quad [g]\mapsto \val(\det(g)).
$$
Let $\xx_{\gamma}^{0}$ be the central connected component of $\xx_{\gamma}$, then $\xx_{\gamma}\cong \xx_{\gamma}^{0}\times \bz$. 
Let $$
\Lambda^{0}=\{\lambda\in \Lambda\mid \val(\det(\lambda))=0\},
$$ 
then $\Lambda^{0}$ acts naturally on $\xx_{\gamma}^{0}$.
The quotient $Z_{\gamma}:=\Lambda^{0}\backslash \xx_{\gamma}^{0}$ is a projective algebraic variety over $k$ and $\xx_{\gamma}^{0}\to Z_{\gamma}$ is an \'etale Galois covering of Galois group $\Lambda^{0}$.

Fix a partition $\{1,\cdots, r\}=J_{1}\sqcup J_{2}$, denoted $J_{\bullet}$, or equivalently fix a character
$$
\kappa:\Lambda^{0}\to \{\pm 1\},\quad \lambda\mapsto (-1)^{\sum_{i\in J_{1}}\lambda_{i}}=(-1)^{\sum_{i\in J_{2}}\lambda_{i}}. 
$$
For $i=1,2$, let $\gl_{J_{i}}=\prod_{j\in J_{i}}\gl_{d_{j}}$, and let $\gamma_{J_{i}}=\diag(\gamma_{j})_{j\in J_{i}}$.
As before, we have the affine Springer fibers $\xx_{\gamma_{J_{i}}}$ and the free discrete abelian groups $\Lambda_{J_{i}}$, and the quotient $Z_{\gamma_{J_{i}}}$.
We have a natural closed embedding 
$$
\xx_{\gamma_{J_{1}}}\times \xx_{\gamma_{J_{2}}}\hookrightarrow \xx_{\gamma}, \quad ([g_{1}], [g_{2}])\mapsto [\diag(g_{1}, g_{2})],
$$
which induces a closed embedding
$$
Z_{\gamma_{J_{1}}}\times Z_{\gamma_{J_{2}}}\hookrightarrow Z_{\gamma}
$$
of codimension
$$
r=\sum_{i\in J_{1}, \,j\in J_{2}} r_{ij},
$$
where $r_{ij}=\dim_{k}(\co[x]/(P_{i}(x), P_{j}(x))).$
As $\xx_{\gamma}^{0}\to Z_{\gamma}$ is an \'etale Galois covering with Galois group $\Lambda^{0}$, the character $\kappa:\Lambda^{0}\to \{\pm 1\}$ defines an $\ell$-adic local system of rank one over $Z_{\gamma}$.

\begin{thm}[Geometric fundamental lemma]\label{main}
There exists a canonical isomorphism
$$
H^{\bullet-2r}(Z_{\gamma_{1}}\times Z_{\gamma_{2}}, \ql)(-r)\cong H^{\bullet}(Z_{\gamma}, \cl). 
$$

\end{thm}

The theorem has been conjectured by Kottwitz. In \cite{laumon unitary}, Laumon gave a proof of this theorem under the purity hypothesis of Goresky, Kottwitz and MacPherson \cite{gkm homology}, which states that $\xx_{\gamma}$ is cohomologically pure in the sense of Grothendieck-Deligne.  
We will give a proof of the theorem based on our work \cite{chen decomposition}. 
As explained in \cite{laumon springer}, \S1.3-1.4, the geometric fundamental lemma implies the arithmetic fundamental lemma for the unitary groups. 
This gives yet another proof of the famous lemma, different from that of Laumon-Ng\^o \cite{laumon ngo}, Ng\^o \cite{ngo}, and Groechenig-Wyss-Ziegler \cite{gwz}.

\section{Proof of the theorem}

We make a brief account of our work \cite{chen decomposition}. According to Laumon \cite{laumon springer}, $Z_{\gamma}$ is homeomorphic to the compactified Jacobian $\overline{P}_{C_{\gamma}}$ of a \emph{spectral curve} $C_{\gamma}$, which is an irreducible projective algebraic curve over $k$, with two points $c$ and $\infty$ such that
\begin{enumerate}[nosep, label=(\roman*)]
\item
$C_\gamma$ has unique singularity at $c$ and $\widehat{\co}_{C_\gamma,\,c}\cong \co[\gamma]$,

\item

the normalization of $C_\gamma$ is isomorphic to $\bp^{1}$ by an isomorphism sending $\infty\in \bp^{1}$ to $\infty\in C_\gamma$.

\end{enumerate}

\noindent 
Let $\pi:(\cc, C_{\gamma})\to (\cb,0)$ be an algebraization of a miniversal deformation of $C_{\gamma}$, and let
$$
f:\overline{\cp}=\overline{\pic}{}^{\,0}_{\cc/ \cb}\to \cb
$$ 
be the relative compactified Jacobian of the family $\pi$.
For any $n\in \bn$, $(n,p)=1$, consider the $\Lambda^{0}/n$-equivariant deformation of the $\Lambda^{0}/n$-covering $\overline{P}_{n}$ of $\overline{P}_{C_{\gamma}}$, or equivalently to extend the $\Lambda^{0}/n$-covering $\overline{P}_{n}\to \overline{P}_{C_{\gamma}}$ over the family $f:\overline{\cp}\to \cb$. 
With the general theory of finite abelian coverings and the auto-duality of compactified Jacobians, we construct the $\Lambda^{0}/n$-coverings 
$$
\Psi_{n}:\cc_{n}\to \cc\times_{\cb}\cb_{n},\quad
\Phi_{n}: \overline{\cp}_{n}\to \overline{\cp}\times_{\cb}\cb_{n},
$$ 
where $\varpi_{n}:\cb_{n}\to  \cb^{\circ}$ is a quasi-finite \'etale covering of an open subscheme $ \cb^{\circ}$ of $\cb$. Compose them with the natural structural morphism to $\cb_{n}$, we get families 
$$
\pi_{n}:\cc_{n}\to \cb_{n}
\quad \text{and}\quad 
f_{n}:\overline{\cp}_{n}\to \cb_{n}.
$$  
Our main theorem in \cite{chen decomposition} concerns the decomposition of the complex $Rf_{n,*}\ql$ into perverse sheaves. To state it, we need to introduce more notations.
For a non-trivial partition $I_{\bullet}=\bigsqcup_{j=1}^{l}I_{j}$ of $\{1,\cdots,r\}$ of length $\ell(I_{\bullet})=l$, we have defined a closed subscheme $\cs_{I_{\bullet}}$ of $\cb^{\circ}$ of codimension 
$$
h_{I_{\bullet}}=\sum_{j\neq j'=1}^{l}\sum_{i\in I_{j}, i'\in I_{j'}}\dim(\co[x]/(P_{i}(x), P_{i'}(x))),
$$
for which the definition is too involved to be recalled here (cf. \cite{chen decomposition}, Def. 4.23). Let $\cs_{I_{\bullet}}^{\circ}$ be the dense open subscheme of $\cs_{I_{\bullet}}$ consisting of the points $s$ over which the fiber $\cc_{s}$ has exactly $h_{I_{\bullet}}$ ordinary double points. Let $\cs_{I_{\bullet, n}}=\cs_{I_{\bullet}}\times_{\cb^{\circ}} \cb_{n}$ and similarly for $\cs_{I_{\bullet}}^{\circ}$. 
The significance of $\cs_{I_{\bullet},n}$ lies in the fact that the union $\bigcup_{I_{\bullet}\vdash \{1,\cdots,r\}}\cs_{I_{\bullet},n}$ is the locus over which the geometric fibers of $f_{n}$ can have multiple irreducible components, and the irreducible components of $\overline{\cp}_{n,s_{n}}$ for a geometric point $s_{n}$ of $\cs^{\circ}_{I_{\bullet}, n}$ can be naturally parametrized by a group $\bbv_{I_{\bullet}}/n$ (cf. \cite{chen decomposition}, theorem 4.27).   

Let $\pi_{I_{\bullet}}^{\circ}: \cc_{I_{\bullet}}^{\circ}:=\cc\times_{\cb}\cs_{I_{\bullet}}^{\circ}\to \cs_{I_{\bullet}}^{\circ}$ be the restriction of the family $\pi:\cc\to \cb$ to $\cs_{I_{\bullet}}^{\circ}$, it is a flat family of projective geometrically irreducible curves with $h_{I_{\bullet}}$ ordinary double points. 
Let $\tilde{\phi}^{\circ}_{I_{\bullet}}:\widetilde{\cc}_{I_{\bullet}}^{\circ}\to \cc_{I_{\bullet}}^{\circ}$ be the normalization, then the family $\tilde{\pi}_{I_{\bullet}}^{\circ}:\widetilde{\cc}_{I_{\bullet}}^{\circ}\to \cs_{I_{\bullet}}^{\circ}$ is a simultaneous resolution of singularities for the family $\pi_{I_{\bullet}}^{\circ}$. 
The pull-back of invertible sheaves along $\tilde{\phi}_{I_{\bullet}}^{\circ}$ defines a morphism
$$
(\tilde{\phi}_{I_{\bullet}}^{\circ})^{*}:\pic_{\cc^{\circ}_{I_{\bullet}}/\cs^{\circ}_{I_{\bullet}}}\to \pic_{\widetilde{\cc}^{\circ}_{I_{\bullet}}/\cs^{\circ}_{I_{\bullet}}},
$$
and this identifies $\pic^{0}_{\widetilde{\cc}^{\circ}_{I_{\bullet}}/\cs^{\circ}_{I_{\bullet}}}$ as the abelian factor of $\pic_{\cc^{\circ}_{I_{\bullet}}/\cs^{\circ}_{I_{\bullet}}}$. 
By construction, the $\Lambda^{0}/n$-covering $\Psi_{n}$ and $\Phi_{n}$ are associated to a constant finite group scheme $\ct_{n}\subset \pic_{{\cc}/\cb}[n]\times_{\cb} \cb_{n}$.
Let $\widetilde{\ct}^{\circ}_{n,I_{\bullet}}\subset \pic_{\widetilde{\cc}^{\circ}_{I_{\bullet}}/\cs^{\circ}_{I_{\bullet}}}[n]\times_{\cs^{\circ}_{I_{\bullet}}}\cs^{\circ}_{I_{\bullet},n}$ be the image under $(\tilde{\phi}_{I_{\bullet}}^{\circ})^{*}$ of the base change to $\cs^{\circ}_{I_{\bullet},n}$ of $\ct_{n}$. The finite group scheme $\widetilde{\ct}^{\circ}_{n, I_{\bullet}}$ defines a finite abelian covering 
\begin{equation*}
\widetilde{\Phi}^{\circ}_{n,I_{\bullet}}:\widetilde{\cp}^{\circ}_{I_{\bullet},n}\to  \pic^{0}_{\widetilde{\cc}^{\circ}_{I_{\bullet}}/\cs^{\circ}_{I_{\bullet}}}\times_{\cs^{\circ}_{I_{\bullet}}}\cs^{\circ}_{I_{\bullet},n}. 
\end{equation*}
Compose it with the structural morphism to $\cs^{\circ}_{I_{\bullet},n}$, we get an Abelian scheme
\begin{equation*}
\tilde{f}^{\circ}_{I_{\bullet},n}: \widetilde{\cp}^{\circ}_{I_{\bullet}, n}\to \cs^{\circ}_{I_{\bullet},n}.
\end{equation*}
Let $\cf_{I_{\bullet},n}^{i}=R^{i}(\tilde{f}^{\circ}_{I_{\bullet},n})_{*}\ql$ and we simplify $\cf_{I_{\bullet}}^{i}=\cf_{I_{\bullet},1}^{i}$. As $\widetilde{\Phi}^{\circ}_{n,I_{\bullet}}$ is an isogeny of abelian schemes, we have actually $\cf_{I_{\bullet},n}^{i}=\varpi_{n}^{*}\cf_{I_{\bullet}}^{i}$. Moreover, we have $\cf_{I_{\bullet}}^{i}=\bigwedge^{i}\cf_{I_{\bullet}}^{1}$ and $\cf_{I_{\bullet}}^{1}\cong R^{1}\tilde{\pi}_{I_{\bullet},*}^{\circ}\ql$.

\begin{thm}[\cite{chen decomposition}, theorem 5.2]\label{support variant}

For the family $f_{n}:\overline{\cp}_{n}\to \cb_{n}$, we have 
$$
Rf_{n,*}\ql=\bigoplus_{i=0}^{2\delta_{\gamma}} j_{n,!*}(R^{\,i}f^{\sm}_{n,*}\ql)[-i]
\oplus 
\bigoplus_{\substack{I_{\bullet} \text{ partition of}\\ \{1,\cdots,r\}}} 
\bigoplus_{i'=0}^{2(\delta_{\gamma}-h_{I_{\bullet}})}
\bigg\{(j_{I_{\bullet},n})_{!*}\cf_{I_{\bullet},n}^{i'}(-h_{I_{\bullet}})[-i'-2h_{I_{\bullet}}]\bigg\}^{\oplus\, |(\bbv_{I_{\bullet}}/n)^{\circ}|},
$$
where $j_{n}:\cb_{n}^{\sm}\to \cb_{n}$ and $j_{I_{\bullet},n}:\cs_{I_{\bullet},n}^{\circ}\to \cs_{I_{\bullet},n}$ are the natural inclusions.

\end{thm}

The first summand in the theorem is the \emph{main term} for the decomposition of $Rf_{n,*}\ql$. 
For $n=1$, the remaining terms vanishes and we recover Ng\^o's support theorem \cite{ngo}.
In general, for the remaining terms, they are of the same nature as the main term. 
For the strata $\cs_{I_{\bullet}}$, let $\varphi_{I_{\bullet}}:\widetilde{\cs}_{I_{\bullet}}\to \cs_{I_{\bullet}}$ be the normalization, we have shown that $\widetilde{\cs}_{I_{\bullet}}$ is smooth over $k$, and there exists a partial resolution $\tilde{\phi}_{I_{\bullet}}^{+}: \widetilde{\cc}^{+}_{I_{\bullet}}\to {\cc}^{+}_{I_{\bullet}}:={\cc}\times_{\cb} \widetilde{\cs}_{I_{\bullet}}$, such that the natural morphism  $\tilde{\pi}^{+}_{I_{\bullet}}: \widetilde{\cc}^{+}_{I_{\bullet}}\to \widetilde{\cs}_{I_{\bullet}}$ is an algebraization of a versal deformation of a partial resolution $\widetilde{C}_{I_{\bullet}}$ of $ C_{\gamma}$. 
Here the partial resolution $\phi_{I_{\bullet}}:\widetilde{C}_{I_{\bullet}}\to C_{\gamma}$ is an isomorphism over $C_{\gamma}\backslash \{c\}$, and such that $\phi_{I_{\bullet}}^{-1}(c)=\{c_{1},\cdots,c_{l}\}$ and that the singularity of $\widetilde{C}_{I_{\bullet}}$ at $c_{j}$ is isomorphic to $\spf(\co[x]/\prod_{i\in I_{j}}P_{i}(x))$ for all $j$.
By construction, the restriction of $\varphi_{I_{\bullet}}$ to $\widetilde{\cs}_{I_{\bullet}}^{\circ}$ is an isomorphism, and with it we can identify  the restriction of the family $\tilde{\pi}_{I_{\bullet}}^{+}$ to $\widetilde{\cs}_{I_{\bullet}}^{\circ}$ with the family $\tilde{\pi}_{I_{\bullet}}^{\circ}$.

Let $\tilde{f}_{I_{\bullet}}^{+}:\widetilde{\cp}_{I_{\bullet}}^{+}\to \widetilde{\cs}_{I_{\bullet}}$ be the relative compactified Jacobian of the family $\tilde{\pi}^{+}_{I_{\bullet}}: \widetilde{\cc}^{+}_{I_{\bullet}}\to \widetilde{\cs}_{I_{\bullet}}$. The partial resolution $\tilde{\phi}_{I_{\bullet}}^{+}: \widetilde{\cc}^{+}_{I_{\bullet}}\to {\cc}^{+}_{I_{\bullet}}$ induces a morphism
$$
(\tilde{\phi}_{I_{\bullet}}^{+})^{*}:\pic_{\cc^{+}_{I_{\bullet}}/\widetilde{\cs}_{I_{\bullet}}}\to \pic_{\widetilde{\cc}^{+}_{I_{\bullet}}/\widetilde{\cs}_{I_{\bullet}}}.
$$ 
Let $\widetilde{\cs}_{I_{\bullet},n}=\widetilde{\cs}_{I_{\bullet}}\times_{\cs_{I_{\bullet}}}\cs_{I_{\bullet}, n}$. 
Note that $\widetilde{\cs}_{I_{\bullet},n}^{\circ}:=\widetilde{\cs}_{I_{\bullet},n}\times_{\cs_{I_{\bullet}, n}}\cs_{I_{\bullet},n}^{\circ}$ coincides with $\cs_{I_{\bullet},n}^{\circ}$ and $\widetilde{\cs}_{I_{\bullet},n}$ is the normalization of ${\cs}_{I_{\bullet},n}$ as $\cs_{I_{\bullet}, n}$ is \'etatle over $\cs_{I_{\bullet}}$.
Let
$$\widetilde{\ct}^{+}_{n,I_{\bullet}}\subset \pic_{\widetilde{\cc}^{+}_{I_{\bullet}}/\widetilde{\cs}_{I_{\bullet}}}[n]\times_{\widetilde{\cs}_{I_{\bullet}}}\widetilde{\cs}_{I_{\bullet},n}
$$
be the image under $(\tilde{\phi}_{I_{\bullet}}^{+})^{*}$ of the base change to $\widetilde{\cs}_{I_{\bullet},n}$ of $\ct_{n}\subset \pic_{{\cc}/\cb}[n]\times_{\cb} \cb_{n}$. 
Then $\widetilde{\ct}^{+}_{n,I_{\bullet}}$ defines a finite abelian covering 
\begin{equation*}
\widetilde{\Phi}^{+}_{n,I_{\bullet}}:\widetilde{\cp}^{+}_{I_{\bullet},n}\to  \widetilde{\cp}^{+}_{I_{\bullet}}\times_{\widetilde{\cs}_{I_{\bullet}}}\widetilde{\cs}_{I_{\bullet},n}. 
\end{equation*}
Compose it with the structural morphism to $\widetilde{\cs}_{I_{\bullet},n}$, we get a family
\begin{equation*}
\tilde{f}^{+}_{I_{\bullet},n}: \widetilde{\cp}^{+}_{I_{\bullet}, n}\to \widetilde{\cs}_{I_{\bullet},n}.
\end{equation*}
By construction, the restriction of $\widetilde{\ct}^{+}_{n,I_{\bullet}}$ to $\widetilde{\cs}_{I_{\bullet},n}^{\circ}$ coincides with $\widetilde{\ct}^{\circ}_{n,I_{\bullet}}$, and the restriction of $\tilde{f}^{+}_{I_{\bullet},n}$ to $\widetilde{\cs}_{I_{\bullet},n}^{\circ}$ coincides with  $\tilde{f}^{\circ}_{I_{\bullet},n}$.

By construction, the curve $\widetilde{C}_{I_{\bullet}}$ is projective rational with singularities at $c_{1},\cdots,c_{l}$, such that $\widehat{\co}_{\widetilde{C}_{I_{\bullet}},c_{j}}\cong \co[\gamma_{I_{j}}]$ for all $j$. 
For $I\subsetneq \{1,\cdots,r\}$, let $C_{\gamma_{I}}$ be the spectral curve of $\gamma_{I}$.  
Let $\pi_{I}:\cc_{I}\to \cb_{I}$ be an algebraization of a miniversal deformation of $C_{\gamma_{I}}$, and let $f_{I}:\overline{\cp}_{I}\to \cb_{I}$ be its relative compactified Jacobian.
Let $j_{I}:\cb_{I}^{\rm sm}\to \cb_{I}$ be the inclusion and let $\cf_{I}^{i}=R^{i}f_{I,*}^{\sm}\ql$. 
We define $f_{I,n}:\overline{\cp}_{I,n}\to \cb_{I,n}$ and denote $j_{I,n}:\cb_{I,n}^{\rm sm}\to \cb_{I,n}$ as before, let $\cf_{I,n}^{i}=R^{i}(f_{I,n}^{\sm})_{*}\ql$. 
Let $M_{I_{\bullet}}$ be the Levi subgroup of $\gl_{d}$ defined by 
\begin{equation*}
M_{I\bullet}=\gl_{I_{1}}\times \cdots \times \gl_{I_{l}}.
\end{equation*}
It is clear that $\gamma$ belongs to the Lie algebra of $M_{I\bullet}$, and we have the affine Springer fiber 
\begin{equation*} 
\xx_{\gamma}^{I_{\bullet}}:=\xx_{\gamma}^{M_{I_{\bullet}}}=\xx_{\gamma_{I_{1}}}\times\cdots\times \xx_{\gamma_{I_{l}}}.
\end{equation*}
With the glueing construction of Laumon \cite{laumon springer}, we have the radical finite surjective morphisms
$$
\Lambda_{I_{j}}^{0}\backslash \xx_{\gamma_{I_{j}}}^{0}\to \overline{P}_{C_{\gamma_{I_{j}}}}, \quad  j=1,\cdots, l,$$
and
\begin{equation*}
\big(\Lambda_{I_{1}}^{0}\backslash \xx_{\gamma_{I_{1}}}^{0}\big)\times\cdots\times \big(\Lambda_{I_{l}}^{0}\backslash \xx_{\gamma_{I_{l}}}^{0}\big)\to \overline{P}_{\widetilde{C}_{I_{\bullet}}}.
\end{equation*}
Hence we get the factorization
\begin{equation}\label{factorize cj si}
\overline{P}_{\widetilde{C}_{I_{\bullet}}} \cong \overline{P}_{C_{\gamma_{I_{1}}}}\times\cdots\times \overline{P}_{C_{\gamma_{I_{l}}}}.
\end{equation}
This factorization property is compatible with the local-global property of the deformation theory of curves and their compactified Jacobians.
In particular, the family $\tilde{f}^{+}_{I_{\bullet}}: \widetilde{\cp}^{+}_{I_{\bullet}}\to \widetilde{\cs}_{I_{\bullet}}$ is isomorphic to the product of versal deformations of $\overline{P}_{C_{\gamma_{I_{j}}}}, j=1,\cdots,l$.
 This implies the factorization
$$
R\tilde{f}^{+}_{I_{\bullet}, *}\ql=(R{f}_{I_{1},*}\ql)\boxtimes\cdots\boxtimes (R{f}_{I_{l},*}\ql).
$$
The factorization (\ref{factorize cj si}), restricted to the Jacobian of the curve, induces a factorization of the finite group scheme $\widetilde{\ct}^{+}_{n,I_{\bullet}}$. Each factor determines a $\Lambda_{I_{j}}^{0}/n$-covering $f_{I_{j},n}: \overline{\cp}_{I_{j}, n}\to \cb_{I_{j},n}$ of the family $f_{I_{j}}: \overline{\cp}_{I_{j}}\to \cb_{I_{j}}$. 
Hence we get the factorization
$$
R(\tilde{f}^{+}_{I_{\bullet}, n})_{*}\ql=\big(R({f}_{I_{1},n})_{*}\ql\big)\boxtimes\cdots\boxtimes \big(R({f}_{I_{l},n}\big)_{*}\ql).
$$
Applying theorem \ref{support variant} to each term at the right, we get:

\begin{thm}\label{reduction theorem}
With the above notations, we have 
\begin{align*}
R(\tilde{f}^{+}_{I_{\bullet}, n})_{*}\ql=\bigboxtimes_{j=1}^{l}&\Bigg[\bigoplus_{i=0}^{2\delta_{\gamma_{I_{j}}}} (j_{I_{j}, n})_{!*}(R^{\,i}(f^{\sm}_{I_{j},n})_{*}\ql)[-i]
\oplus 
\bigoplus_{I_{j, \bullet} \text{ partition of }I_{j}} 
\\&
\bigoplus_{i'=0}^{2(\delta_{\gamma_{I_{j}}}-h_{I_{j,\bullet}})}
\bigg\{(j_{I_{j,\bullet},n})_{!*}\cf_{I_{j,\bullet},n}^{i'}(-h_{I_{j,\bullet}})[-i'-2h_{I_{j,\bullet}}]\bigg\}^{\oplus\, |(\bbv_{I_{j,\bullet}}/n)^{\circ}|}\Bigg],
\end{align*}
where the addition of $I_{j}$ at the subscript means similar objects for the family $f_{I_{j},n}: \overline{\cp}_{I_{j}, n}\to \cb_{I_{j},n}$.

\end{thm}

As explained before proposition 5.4 of \cite{chen decomposition} and during its proof, the group $\Lambda^{0}/n$ acts on each summand $((j_{I_{\bullet}, n})_{!*}\cf_{I_{\bullet},n}^{i})^{(\bbv_{I_{\bullet}}/n)^{\circ}}$ in theorem \ref{support variant}, and the action is induced by the action on the set of irreducible components of the geometric fibers over $\cs_{I_{\bullet}, n}^{\circ}$, which equals $\bbv_{I_{\bullet}}/n$ by theorem 4.27 of \cite{chen decomposition}.  
We can make the action more explicit. For simplicity, we denote $\Lambda_{I_{\bullet}}^{0}=\Lambda_{I_{1}}^{0}\times\cdots\times\Lambda_{I_{l}}^{0}$.

\begin{lem}\label{eigen decomposition}

We have canonical isomorphism $\bbv_{I_{\bullet}}/n\cong (\Lambda^{0}/n)\big/(\Lambda_{I_{\bullet}}^{0}/n)$. 

\end{lem}

\begin{proof}

Let ${\varphi}_{I_{\bullet},n}: \widetilde{\cs}_{I_{\bullet},n}\to \cs_{I_{\bullet}, n}$ be the structural morphism, it is a finite morphism as it is the base change of the normalization $\varphi_{I_{\bullet}}:\widetilde{\cs}_{I_{\bullet}}\to  \cs_{I_{\bullet}}$. 
Combined with the fact that the restriction of $\tilde{f}^{+}_{I_{\bullet},n}$ to $\widetilde{\cs}_{I_{\bullet},n}^{\circ}$ coincides with  $\tilde{f}^{\circ}_{I_{\bullet},n}$, we get
\begin{align}
(j_{I_{\bullet}, n})_{!*}\cf_{I_{\bullet},n}^{i}&=\big({\varphi}_{I_{\bullet},n}\big)_{*} 
\big(\tilde{j}^{+}_{I_{\bullet},n}\big)_{!*}R^{\,i}\big(\tilde{f}^{+, \sm}_{I_{\bullet},n}\big)_{*}\ql \nonumber
\\
&=\big({\varphi}_{I_{\bullet},n}\big)_{*}\bigboxtimes_{i_{1}+\cdots+i_{l}=i} (j_{I_{j}, n})_{!*}(R^{\,i_{j}}(f^{\sm}_{I_{j},n})_{*}\ql). \label{reduce 2}
\end{align}
Apply proposition 5.4 of \cite{chen decomposition} to the family $f_{I_{j},n}$, we obtain that $\Lambda_{I_{j}}^{0}/n$ acts trivially on $(j_{I_{j}, n})_{!*}R(f^{\sm}_{I_{j},n})_{*}\ql$.
With the isomorphism (\ref{reduce 2}), we obtain that the action of $\Lambda^{0}/n$ on $((j_{I_{\bullet}, n})_{!*}\cf_{I_{\bullet},n}^{i})^{(\bbv_{I_{\bullet}}/n)^{\circ}}$ factors through $(\Lambda^{0}/n)\big/(\Lambda_{I_{\bullet}}^{0}/n)$.
On the other hand, the action of $\Lambda^{0}/n$ on the irreducible components of the geometric fibers over $\cs_{I_{\bullet}, n}^{\circ}$ must be transitive. Indeed, $\Lambda^{0}/n$ acts transitively on the set of irreducible components of $\overline{P}_{n}$, hence also on $(R^{2\delta_{\gamma}}f_{n,*}\ql)_{0_{n}}$, here $0_{n}\in \cb_{n}$ is the point lying under $\overline{P}_{n}$. Now it suffices to invoke theorem \ref{support variant} to arrive at the conclusion. 

With the above results, we get a surjective homomorphism of groups
$$
(\Lambda^{0}/n)\big/(\Lambda_{I_{\bullet}}^{0}/n)\to \bbv_{I_{\bullet}}/n.
$$
By construction, $\bbv_{I_{\bullet}}$ is the free abelian group generated by the classes $c_{I_{j}}$ which satisfies $c_{I_{1}}+\cdots+c_{I_{l}}=0$ (cf. \cite{chen decomposition}, \S4.3).  
Hence the groups on the source and target of the above morphism have the same cardinal, they have to be isomorphic.

\end{proof}

We can now proceed to the proof of theorem \ref{main}.
To calculate $H^{*}(Z_{\gamma}, \cl)=H^{*}(\overline{P}_{C_{\gamma}}, \cl)$, we use the Galois covering $\overline{P}_{n}\to \overline{P}_{C_{\gamma}}$. Here we are taking $\ell$ sufficiently large and $n<\ell$ such that $n$ is even and coprime to both $p$ and $\ell$. 
The character $\kappa:\Lambda^{0}\to \{\pm 1\}$ factors through $\Lambda^{0}/n$ and we get a character of the latter, denoted $\kappa_{n}$. Then the local system is the push-forward of the covering $\overline{P}_{n}\to \overline{P}_{C_{\gamma}}$ via the character $\kappa_{n}$, and so we have a spectral sequence
$$
E_{2}^{p,q}=H^{p}(\Lambda^{0}/n, H^{q}(\overline{P}_{n}, \ql))\Longrightarrow H^{p+q}(\overline{P}_{C_{\gamma}}, \cl).
$$ 
As the Galois group $\Lambda^{0}/n$ is finite, $E_{2}^{p,q}=0$ for $p\neq 0$, the spectral sequence degenerates and we get
\begin{equation}\label{reduce to kappa 0}
H^{*}(\overline{P}_{C_{\gamma}}, \cl)=H^{*}(\overline{P}_{n}, \ql)_{\kappa_{n}},
\end{equation}
where the subscript ${}_{\kappa_{n}}$ denotes the $\kappa_{n}$-isotypical subspace.
Applying theorem \ref{support variant} to the right hand side of the equation (\ref{reduce to kappa 0}), we get
\begin{align}
H^{i}(\overline{P}_{C_{\gamma}}, \cl)=&\bigoplus_{i'=0}^{i} \ch^{i-i'}\big(j_{n,!*}(R^{\,i'}f^{\sm}_{n,*}\ql)\big)_{0_{n},\kappa_{n}}
\oplus \nonumber
\\
&\bigoplus_{\substack{I_{\bullet} \text{ partition of}\\ \{1,\cdots,r\}}} \bigoplus_{i''=0}^{i-2h_{I_{\bullet}}}
\bigg\{\ch^{i-2h_{I_{\bullet}}-i''}\big((j_{I_{\bullet},n})_{!*}\cf_{I_{\bullet},n}^{i''}\big)(-h_{I_{\bullet}})^{\oplus\, |(\bbv_{I_{\bullet}}/n)^{\circ}|}\bigg\}_{0_{n}, \kappa_{n}}.  \label{reduce to kappa}\end{align}
We can analyse the right hand side of the equation with lemma \ref{eigen decomposition}. 
By the lemma, the action of $\Lambda^{0}/n$ on $\bbv_{I_{\bullet}}/n$ factors through $(\Lambda^{0}/n)\big/(\Lambda_{I_{\bullet}}^{0}/n)$, and is essentially the regular representation. In particular, the group $\Lambda_{I\bullet}^{0}/n$ acts trivially on the summand indexed by $I_{\bullet}$. Compare it with $\Ker(\kappa_{n})$. 
For the partition $I_{\bullet}$ which is not a refinement of the partition $J_{\bullet}$, there is an $j_{0}$ such that $I_{j_{0}}\cap J_{1}\neq \emptyset$ and $I_{j_{0}}\cap J_{2}\neq \emptyset$, from which we can find an element $(\cdots, 1,\cdots, -1,\cdots)\in (\Lambda_{I_{1}}^{0}/n)\times\cdots\times (\Lambda_{I_{l}}^{0}/n)$ but not in $\Ker(\kappa_{n})$. This implies that the contributions of the summands indexed by such partitions are $0$.
Similarly, for the partition $I_{\bullet}$ which is a strict refinement of $J_{\bullet}$, we have $\Lambda^{0}_{I_{\bullet}}\subsetneq \Lambda_{J_{\bullet}}^{0}\subset\Ker(\kappa_{n})$, whence the contribution of such summands are $0$ as well. 
For the summand indexed by $J_{\bullet}$, recall that $(\bbv_{I_{\bullet}}/n)^{\circ}=(\bbv_{I_{\bullet}}/n)\backslash \bigcup_{I_{\bullet}\vdash I_{\bullet}'}\bbv_{I_{\bullet}'}/n$, where $I_{\bullet}\vdash I_{\bullet}'$ means that $I_{\bullet}$ is a strict  refinement of $I_{\bullet}'$, we obtained that $(\bbv_{J_{\bullet}}/n)^{\circ}$ as a representation of $(\Lambda^{0}/n)\big/(\Lambda_{J_{\bullet}}^{0}/n)$ consists of all the non-trivial representations in the regular representation. 
Hence the summand indexed by $J_{\bullet}$ in the equation (\ref{reduce to kappa}) is exactly of dimension $1$, and we conclude
\begin{equation}\label{final 0}
H^{*}(\overline{P}_{C_{\gamma}}, \cl)=\bigoplus_{i=0}^{2(\delta_{\gamma}-h_{J_{\bullet}})}
\Big\{(j_{J_{\bullet},n})_{!*}\cf_{J_{\bullet},n}^{i}(-h_{J_{\bullet}})[-i-2h_{J_{\bullet}}]\Big\}_{0_{n}}\otimes V_{\kappa_{n}}.
\end{equation}
Note that $h_{J_{\bullet}}=r$ by definition, and $\otimes V_{\kappa_{n}}$ doesn't change anything if we consider the equation as an equality of abstract groups. 
Moreover, as we have mentioned, $\cf_{J_{\bullet},n}^{i}=\varpi_{n}^{*}\cf_{J_{\bullet}}^{i}$ for the \'etale morphism $\varpi_{n}:\cb_{n}\to \cb^{\circ}$, and $j_{J_{\bullet},n}:\cb_{n}^{\sm}\to \cb_{n}$ is the base change of $j_{J_{\bullet}}:\cb^{\sm}\to \cb$ via $\varpi_{n}$. Hence the equation (\ref{final 0}) simplifies to  
\begin{align*}
H^{*}(\overline{P}_{C_{\gamma}}, \cl)=\bigoplus_{i=0}^{2(\delta_{\gamma}-r)}
\Big\{(j_{J_{\bullet}})_{!*}\cf_{J_{\bullet}}^{i}(-r)[-i-2r]\Big\}_{0}
=H^{*-2r}(\overline{P}_{C_{\gamma_{J_{1}}}}\times \overline{P}_{C_{\gamma_{J_{2}}}}, \ql)(-r).
\end{align*}
Here for the second equality we are using the equation (\ref{reduce 2}) for the case $n=1$ and the support theorem of Ng\^o \cite{ngo} for $\overline{P}_{C_{\gamma_{J_{i}}}}, i=1,2$.
As $\overline{P}_{C_{\gamma}}$ is homeomorphic to $Z_{\gamma}$, this finishes the proof of theorem \ref{main}.

\bigskip
\small
\noindent
\begin{tabular}{ll}
&Zongbin {\sc Chen} \\ 
\\
&School of mathematics, Shandong University\\
&250100, JiNan, Shandong, \\
&P. R. China \\
&email: {\tt zongbin.chen@email.sdu.edu.cn}

\end{tabular}


\begin{thebibliography}{100}
\labelwidth=4em
\addtolength\leftskip{25pt}
\setlength\labelsep{0pt}
\addtolength\parskip{\smallskipamount}






\bibitem[C]{chen decomposition}{Z. Chen, \textit{A decomposition theorem for the affine Springer fibers}. \url{https://arxiv.org/pdf/2404.08225}.}


\bibitem[GKM]{gkm homology}{M. Goresky, R. Kottwitz, R. MacPherson, \textit{Homology of affine Springer fibers in the unramified case}, Duke Math. J. 121 (2004), no. 3, 509-561.}

\bibitem[GWZ]{gwz}{M. Groechenig, D. Wyss, P. Ziegler, 
\textit{Geometric stabilisation via p-adic integration}.
J. Amer. Math. Soc.33(2020), no.3, 807-873.}


\bibitem[L1]{laumon springer}{G. Laumon, \textit{Fibres de Springer et jacobiennes compactifi\'ees}, \url{https://arxiv.org/pdf/math/0204109}.}

\bibitem[L2]{laumon unitary}{G. Laumon, \textit{Sur le lemme fondamental pour les groupes unitaires}, \url{https://arxiv.org/pdf/math/0212245}}



\bibitem[LN]{laumon ngo}{G. Laumon, B. C. Ng\^{o}, \textit{Le lemme fondamental pour les groupes unitaires}, Ann. Math., 168 (2008), 477-573.}

\bibitem[N]{ngo}{Bau Ch\^{a}u Ng\^o, \textit{Le lemme fondamental pour les alg\`ebres de Lie.} Publ. Math. IHES. No. 111 (2010), 1–169.}


\end{thebibliography}
\end{document}